\documentclass[12pt]{amsart}
\usepackage{amssymb}

\usepackage[cmtip,matrix,arrow]{xy}

\newtheorem{theorem}{Theorem}[section]
\newtheorem{corollary}[theorem]{Corollary}
\newtheorem{proposition}[theorem]{Proposition}

\newtheorem{lemma}[theorem]{Lemma}

\theoremstyle{remark}
\newtheorem{remark}[theorem]{Remark}

\newtheorem{remarks}[theorem]{Remarks}
\newtheorem{example}[theorem]{Example}

\newcommand\A{\mathcal{A}}
\newcommand\be{\begin{equation}}
\newcommand\ee{\end{equation}}
\newcommand\M{\mathcal{M}}
\renewcommand\L{\mathcal{L}}
\newcommand\G{\mathcal{G}}

\newcommand{\Co}{\mathcal{C}}

\newcommand{\U}{\on{U}}

\newcommand{\R}{\mathbb{R}}
\newcommand{\C}{\mathbb{C}}

\newcommand{\Z}{\mathbb{Z}}

\newcommand\lie[1]{\mathfrak{#1}}
\renewcommand{\k}{\lie{k}}

\newcommand{\g}{\lie{g}}
\newcommand{\m}{\lie{m}}

\renewcommand{\t}{\lie{t}}
\newcommand{\Alc}{\lie{A}}

\newcommand{\on}{\operatorname}

\newcommand{\Ad}{ \on{Ad} } \newcommand{\Hol}{ \on{Hol} }

\newcommand{\Hom}{ \on{Hom}} 
 
 \newcommand{\Spin}{ \on{Spin}}
\newcommand{\SU}{ \on{SU}} 
\newcommand{\SO}{ \on{SO}}

\newcommand{\Mult}{ \on{Mult}} \newcommand{\Vol}{ \on{Vol}}
 
\newcommand{\Oma}{{\underline \Om}}

\newcommand{\D}{ \mathcal{D} }

\newcommand\dirac{/\kern-1.2ex\partial} 
\newcommand\qu{/\kern-.7ex/} 
\newcommand{\fus}{\circledast}  

\newcommand{\labell}\label

\newcommand{\rb}[1]{\raisebox{1.5ex}[-1.5ex]{#1}}

\newcommand{\hra}{\hookrightarrow}

\renewcommand{\d}{{\mbox{d}}}
\newcommand{\ol}{\overline}

\newcommand{\pr}{\on{pr}}
\newcommand\lam{\lambda}

\newcommand\Sig{\Sigma}

\newcommand\eps{\epsilon}
\newcommand\Om{\Omega}
\newcommand\om{\omega}

\newcommand{\f}{\frac}

\newcommand{\p}{\partial}

\newcommand{\hh}{{\textstyle \f{1}{2}}}

\newcommand{\olt}{\overline{\theta}}

\newcommand\LG{\widehat{LG}}      
\newcommand\beqn{\begin{equation}}      
\newcommand\eeqn{\end{equation}}      

\newcommand{\wh}{\widehat}

\newcommand{\mf}{\mathfrak}
\newcommand{\beq}{\begin{eqnarray*}}
\newcommand{\eeq}{\end{eqnarray*}}
\include{verlinde.aux}
\begin{document}

\title[Formulas of Verlinde type for non-simply connected
groups]{Formulas of Verlinde type \\ for non-simply connected groups}

\date{\today}

\author{A. Alekseev}
\address{Institute for Theoretical Physics \\ Uppsala University \\
Box 803 \\ \mbox{S-75108} Uppsala \\ Sweden}
\email{alekseev@teorfys.uu.se}

\author{E. Meinrenken}
\address{University of Toronto, Department of Mathematics,
100 St George Street, Toronto, Ontario M5R3G3, Canada }
\email{mein@math.toronto.edu}

\author{C. Woodward}
\address{Mathematics-Hill Center, Rutgers University,
110 Frelinghuysen Road, Piscataway NJ 08854-8019, USA}
\email{ctw@math.rutgers.edu}

\begin{abstract}
We derive Verlinde's formula from the fixed point formula
for loop groups proved in the companion paper \cite{al:fi},
and extend it to compact, connected groups that are not necessarily
simply-connected.
\end{abstract}

\maketitle



\section{Introduction}

In this paper we give applications of the fixed point formula proved
in the companion paper \cite{al:fi}.  Our original motivation was to
understand a formula of E. Verlinde \cite{ve:fr} for the geometric
quantization of the moduli space of flat connections on a Riemann
surface.  In particular A. Szenes suggested to us that the Verlinde
formula should follow from an equivariant index theorem, much as the
Weyl or Steinberg formulas can be interpreted as fixed point formulas
for flag varieties.  It turns out that our approach also applies to
moduli spaces for compact, connected groups that are not necessarily
simply connected.
In the present paper we will consider the case of at most one marking.
The main result is Theorem \ref{th:ns} below.  The case $G=SO(3)$ is
due to Pantev \cite{pa:cm}, and $G = \on{PSU}(p)$ with $p$ prime to
Beauville \cite{be:pg}. A formula for any number of markings was
conjectured very recently by Fuchs and Schweigert \cite{fu:ac}.  The
proof of this more general result, which involves so-called orbit Lie
algebras, will appear in a later work.

Verlinde's formula appears in the literature in various guises.  Our
version computes the index of the pre-quantum line bundle over the
moduli space of flat connections.  Verlinde's original conjecture
computes the dimension of the space of conformal blocks, for a wide
range of two-dimensional conformal field theories.  The cases we
consider here arise from the Wess-Zumino-Witten model; see Felder,
Gawedzki, and Kupiainen \cite{fe:sp} for arbitrary simple groups.  In
the algebraic geometry literature, many authors refer to Verlinde's
formula as computing the dimension of the space of holomorphic
sections (generalized theta-functions).  
In many cases, the higher cohomology of the pre-quantum
line bundle is known to vanish (see recent work by Teleman
\cite{te:qu}) so that our version and the algebraic geometry version
are the same.

Mathematically rigorous approaches to Verlinde's formula in algebraic
geometry are due to Tsuchiya-Ueno-Yamada \cite{ts:cf}, Faltings
\cite{fa:vl}, and Teleman \cite{te:lc,te:ve}, to name a few.  The
comparison between conformal blocks and holomorphic sections can be
found in Kumar-Narasimhan-Ramanathan \cite{ku:in}, Beauville-Laszlo
\cite{be:cb}, and Pauly \cite{pa:co}. A nice survey can be found in
Sorger \cite{so:ve}.

A proof of the Verlinde formula via the Riemann-Roch theorem was
outlined by Szenes in \cite{sz:co}, and carried out for $\SU(2)$ and
$\SU(3)$.  The proof was extended by Jeffrey-Kirwan \cite{jk:in} to
$\SU(n)$, and Bismut-Labourie \cite{bi:sy} to arbitrary compact,
connected simply-connected groups, for sufficiently high level.
The idea of deriving the
Verlinde formula from localization also appears in the physics papers
by Gerasimov \cite{ge:ve} and Blau and Thompson \cite{bl:eq}.

\vskip.3in

\begin{center}
{\bf Notation.}
\end{center}
\vskip.1in
This paper is designed to be read parallel to its companion {\it A
fixed point formula for loop group actions} \cite{al:fi}.  Throughout
this paper, $G$ will denote a compact, connected, simply connected Lie
group, and $G=G_1\times\cdots\times G_s$ its decomposition into simple
factors. Given a tuple $k=(k_1,\ldots,k_s)$ of positive numbers,
$B=B_k$ will denote the invariant bilinear form which restricts to
$k_j$ times the basic inner product on the $j$th factor (\cite{al:fi}, Section
2.4).

We use the following Sobolev norms.  Fix a number $f > 1$.  For any
manifold $X$ (possibly with boundary) and $p\le \dim X$, we denote by
$\Oma^p(X,\g)$ the set of $\g$-valued $p$-forms of Sobolev class
$f-p+\dim X/2$.  Under this assumption, forms in $\Oma^0(X,\g)$ are
$C^1$ and those in $\Oma^1(X,\g)$ are $C^0$.  The space $\Oma^0(X,\g)$
is the Lie algebra of the group $\G(X)=\on{Map}(X,G)$ of Sobolev class
$f+\dim X/2$.  

Some other notations introduced in this paper are
\begin{tabbing}
\indent \= $M_\mu$, $M_0=M\qu G$\quad \= \kill 
\> $\M(\Sigma)$ \> moduli space of flat $G$-connections on 
   surface $\Sig$; \ref{sec:mod} \\
\> $M(\Sigma)$ \> holonomy manifold of $\M(\Sigma)$; \ref{sec:hol} \\
\> $L(\Sigma)$ \> pre-quantum line bundle over $\M(\Sigma)$;
   \ref{sec:pre} \\ 
\> $h$, $r$ \> genus, number of boundary components
   of $\Sigma$; \ref{sec:mod}\\
\>  $\fus$    \>  fusion product; Appendix \ref{sec:fusion}   \\
\>  $\Gamma$, $G'$    \> central subgroup of $G$, quotient $G'=G/\Gamma$; 
    \ref{sec:prime}    \\
\>  $F_\gamma$    \> fixed point component corresponding to $\gamma
    \in \Gamma^{2h}$; \ref{sec:fixed}    \\
\>      \>     \\
\end{tabbing}

\section{The simply-connected case}
\subsection{The moduli space of flat connections}
\label{sec:mod}
We begin with a brief review of the gauge theory construction of
moduli spaces of flat connections. More details can be found in
\cite{at:ge}, \cite{me:lo}, and \cite{do:bv}.  Let $\Sig=\Sig_h^r$
denote a compact, connected, oriented surface of genus $h$ with $r$
boundary components. Let $\A(\Sig)=\Oma^1(\Sig,\g)$ be the affine
space of connections on the trivial $G$-bundle over $\Sig$, equipped
with the action of $\G(\Sig)$ by gauge transformations
\begin{equation}
g\cdot A=\Ad_g(A)-g^*\olt,
\end{equation}
where $\olt$ is the right-invariant Maurer-Cartan form.  Let
$\G_\p(\Sig)\subset \G(\Sig)$ be the kernel of the restriction map
$\G(\Sig)\to \G(\p\Sig)$. Since $G$ is simply connected, the
restriction map is surjective, and therefore
$\G(\Sig)/\G_\p(\Sig)\cong \G(\p\Sig)$.  We define
\begin{equation}\label{eq:moduli}
\M(\Sig):=\A_{\on{flat}}(\Sig)/\G_\p(\Sig),
\end{equation}
the moduli space of flat $G$-connections under based gauge
equivalence.  If $\p\Sig\not=\emptyset$, it is a 
smooth $\G(\p\Sig)$-equivariant Banach manifold.
Pull-back of connections to the boundary induces a map,
$$ \wh{\Phi}:\,\M(\Sig)\to \Oma^1(\p\Sig,\g).
$$
The map $\wh{\Phi}$ is smooth and proper, and is equivariant for the
gauge action of $\G(\p\Sig)$.
Let $B=B_k$ be an invariant inner product on $\g$.  The symplectic
form on $\A(\Sig)$ is given by the integration pairing of 1-forms $
(a_1,a_2) \mapsto \int_\Sigma B(a_1 \wedge a_2).$ As observed by
Atiyah-Bott \cite{at:ge}, the action of $\G_\p(\Sig)$ is Hamiltonian,
with moment map the curvature.  Hence \eqref{eq:moduli} is a
symplectic quotient and $\M(\Sig)$ inherits a symplectic 2-form
$\wh{\om}$.  Moreover, the residual action of $\G(\p\Sig)$ on $\M(\Sig)$ is
Hamiltonian with moment map $\wh{\Phi}$, using the pairing of 
$\Oma^1(\p\Sig,\g)$ and $\Oma^0(\p\Sig,\g)$ given by the inner product 
and integration over $\p\Sig$. 
A choice of parametrization of the boundary $\p \Sig=(S^1)^r$ induces
isomorphisms
$$ \G(\p\Sig) \cong LG^r, \ \ \ \Oma^1(\p\Sig,\g) \cong
\Oma^1(S^1,\g^r) .$$
Thus
$(\M(\Sig),\wh{\om},\wh{\Phi})$ is an example of a Hamiltonian
$LG^r$-manifold with proper moment map.  For any
$\mu=(\mu_1,\ldots,\mu_r)\in L(\g^r)^*$, the symplectic quotient
$\M(\Sig)_\mu$ is the moduli space of flat connections for which the
holonomy around the $j$th boundary component is contained in the
conjugacy class of $\Hol(\mu_j)$. This also covers the case without
boundary, since $M(\Sig_h^0)=M(\Sig_h^1)_0$.

Occasionally we will also use the notation $\M(\Sig,G)$, in order to
indicate the structure group. The decomposition into simple factors
$G=G_1\times \cdots \times G_s$ defines a decomposition of the loop
group $LG=LG_1 \times \cdots \times LG_s$, and the moduli space is the
direct product
$$\M(\Sig,G)=\M(\Sig,G_1)\times \cdots \times \M(\Sig,G_s).$$

\subsection{Pre-quantization of the moduli space}
\label{sec:pre}
The space $\M(\Sig)$ is pre-quantizable at integer level, that is if
all $k_i$ are integers (see e.g. Section 3.3. of \cite{me:lo} or
\cite{fr:cl}).  For later use, we recall the construction of the
pre-quantum line bundle.  The central extension $\wh{\G}(\Sig)$ of
$\G(\Sig)$ is defined by the cocycle
\begin{equation}\label{eq:cocycle}
c(g_1,g_2)=\exp(i\pi \int_\Sig B(g_1^*\theta,g_2^*\olt)).
\end{equation}
The group $\wh{\G}(\Sig)$ acts on the trivial line bundle over $\A(\Sig)$ by 
\begin{equation}\label{eq:lineaction}
(g,z)\cdot (A,w)= \big(g\cdot A, \exp(-i\pi \int_\Sig
B(g^*\theta,A))zw\big).
\end{equation}
The 1-form $a\mapsto \hh \int_\Sig B(A,a)$ on $\A(\Sig)$ defines an
invariant pre-quantum connection.  A trivialization of $\wh{\G}(\Sig)$
over the subgroup $ \G_\p(\Sig)$ is given by the map
\begin{equation}\label{eq:trivialization} 
\alpha:\,\G_\p(\Sig)\to \U(1),\ \ \alpha(g)=\exp(2\pi i
\int_{\Sig\times [0,1]} \ol{g}^*\eta).
\end{equation}
Here $\eta$ is the canonical 3-form on $G$, and
$\ol{g}\in\G(\Sig\times [0,1])$ is any extension such that $\ol{g}=g$
on $\Sig\times\{0\}$ and $\ol{g}=e$ on $(\Sig\times\{1\})\cup
(\p\Sig\times [0,1])$. The map $\alpha$ is well-defined and satisfies
the coboundary condition
$\alpha(g_1g_2)=\alpha(g_1)\alpha(g_2)c(g_1,g_2)$.  One defines the
pre-quantum line bundle as a quotient
$L(\Sig)=(\A_{\on{flat}}(\Sig)\times\C)/\G_\p(\Sig)$; it comes equipped
with an action of $\wh{LG}=\wh{\G}(\Sig)/\G_\p(\Sig)$.

We are thus in the setting of the fixed point formula (\cite{al:fi},
Theorem 4.3) which gives a formula for the $\Spin_c$-index
$\chi(\M(\Sig)_\mu)$.  To apply the fixed point formula we have to (i)
describe the holonomy manifold $M(\Sig):=\M(\Sig)/\Om G^r$, (ii)
determine the fixed point manifolds for elements
$(t_{\lambda_1},\ldots,t_{\lambda_r})$, and (iii) evaluate the fixed
point data. These steps will be carried out in the subsequent
sections.

\subsection{Holonomy manifolds}
\label{sec:hol}
The holonomy manifold $M(\Sig):=\M(\Sig)/\Om G^r$ can be interpreted
as the moduli space of flat connections
$$ M(\Sig)=\A_{\on{flat}}(\Sig)/\{g\in\G(\Sig)|\,g(p_1)=\ldots=g(p_r)=e\}$$
where $p_1,\ldots,p_r$ are the base points on the boundary circles.
The group-valued moment map $\Phi:\,M(\Sig)\to G^r$ takes an
equivalence class of flat connections to its holonomies around the
boundary circles. The 2-form $\om$ has the following explicit
description (see \cite[Section 9]{al:mom}.)  We begin with the case of
a 2-holed sphere $\Sig_0^2$.  The surface $\Sig_0^2$ is obtained from
a 4-gon by identifying the sides according to the word
$D_1AD_2A^{-1}$. Parallel transport along the paths $A$ and
$A^{-1}D_1$ defines a diffeomorphism
$$M(\Sig_0^2)=G\times G.$$
The $G^2$-action is given by
\begin{equation}\label{eq:action} 
(g_1,g_2)\cdot(a,b)=(g_1 a g_2^{-1},\,g_2 b g_1^{-1}).
\end{equation}
The moment map is 
\begin{equation}\label{eq:moment} 
\Phi(a,b)=(ab,a^{-1}b^{-1})
\end{equation}
and the 2-form is given by
\begin{equation}\label{eq:2form}
\om=\hh\big(B(a^*\theta,b^*\olt)+B(a^*\olt,b^*\theta)\big).
\end{equation}

The holonomy manifolds for the general case $\Sig=\Sig_h^r$
are obtained from $M(\Sig_0^2)$ by {\it fusion}, which 
we recall in Appendix \ref{sec:fusion}. First, the moduli space 
$M(\Sig_1^1)$ for the 1-punctured torus is
$$ M(\Sig_1^1)=M(\Sig_0^2)_{\on{fus}}\cong G^2.$$
The $G$-action is conjugation on each factor and the moment map is the
Lie group commutator $\Phi(a,b)=[a,b]=aba^{-1}b^{-1}$.  The moduli
space for the surface of genus $h$ with $1$ boundary component is an
$h$-fold fusion product
$$ M(\Sig_h^1)=M(\Sig_1^1)\fus\ldots\fus M(\Sig_1^1) = G^{2h}.$$
$G$ acts by conjugation on each factor, and the moment map is a
product of Lie group commutators. The moduli space for the $r$-holed
sphere $\Sig_0^r$ is an $(r-1)$-fold fusion product
$$ M(\Sig_0^r)=M(\Sig_0^2)\fus\ldots\fus
M(\Sig_0^2) = G^{2(r-1)}$$
where we fuse with respect to the first $G$-factor for each
$G^2$-space $M(\Sig_0^2)$.  Finally, the moduli space for $\Sig_h^r$
is
$$ M(\Sig_h^r)=M(\Sig_h^1)\fus M(\Sig_0^r) = G^{2(h+r-1)} .$$ 

\subsection{The fixed point sets}

The fixed point sets for the action on the holonomy manifold are
symplectic tori:

\begin{proposition}  \label{prop:fpm}
The fixed point set for the action of 
$(t_{\lam_1},\ldots,t_{\lam_r})$ on $M(\Sig_h^r) = G^{2(h+r-1)}$ is
empty unless $\lambda_1=\ldots=\lambda_r=:\lambda$, and
$$M(\Sig_h^r)^{(t_\lam,\ldots,t_\lam)}=F:=T^{2(h+r-1)}.$$
\end{proposition}

\begin{proof} 
Since $M(\Sig_h^r)$ is obtained from a direct product 
of $h+r-1$ copies of $M(\Sig_0^2)$ by passing to diagonal actions 
for some of the $G$-factors, it suffices to prove 
 Proposition \ref{prop:fpm} for $\Sig_0^2$. By \eqref{eq:action}, an
element $(a,b)\in M(\Sig_0^2)$ is fixed by $(t_{\lam_1},t_{\lam_2})$
if and only if
\begin{equation} \label{eq:teq}
 t_{\lam_1}=\Ad_a t_{\lam_2},\ \ t_{\lam_2}=\Ad_b t_{\lam_1}.
\end{equation}
Both $t_{\lam_1}$ and $t_{\lambda_2}$ belong to the exponential of the
alcove $\exp(\Alc)$.  Since each conjugacy class meet $\exp(\Alc)$
only once, \eqref{eq:teq} holds if and only if $\lambda_1 =
\lambda_2$.
\end{proof}

Notice that the fixed point set is independent of $\lambda$; in fact,
$F$ is fixed by the full diagonal torus $T\subset G^r$.  

\subsection{Evaluation of the fixed point contributions}

Let $\Sig=\Sig_h^r$ and $\mu=(\mu_1,\ldots,\mu_r)$ with
$\mu_j\in\Lambda^*_k$.  By Theorem 4.3. of \cite{al:fi}, the $\Spin_c$-index
is given by the formula
\begin{equation}\label{eq:fp}
\chi(\M(\Sig)_\mu)
=\f{1}{(\# T_{k+c})^r}
\sum_{\lambda\in\Lambda^*_k}
|J(t_\lam)|^{2r} \prod_{j=1}^r \chi_{\mu_j}(t_\lam)^*
\zeta_F(t_\lam)^{1/2} \int_F 
\f{\hat{A}(F)e^{\hh c_1(\L_F)}}{\D_\R(\nu_F,t_\lam)}.
\end{equation}
Here we abbreviated $(t_\lam,\ldots,t_\lam)$ to $t_\lam$, viewing $T$
as diagonally embedded into $G^r$. Since the normal
bundle $\nu_F$ is $T$-equivariantly isomorphic to $(\g/\t)^{2(h+r-1)}$, 
we have
\begin{equation}\label{eq:d}
D_\R(\nu_F,t_\lam)=J(t_\lam)^{2(h+r-1)} 
=(-1)^{(h+r-1)\#\mf{R}} |J(t_\lam)|^{2(h+r-1)}.
\end{equation}
Furthermore, since $F$ is a product of tori we have 
\begin{equation}\label{eq:a}
\hat{A}(F)=1.
\end{equation}
It remains to work out the integral $\int_F \exp(\hh c_1(\L_F))$
and to calculate the phase factor 
$\zeta_F(t_\lam)^{1/2}$.

\begin{proposition}  \label{prop:int} The integral of $ \exp(\hh c_1(\L_F))$
over $F$ equals ${(\# T_{k+c})}^{h+r-1}.$
\end{proposition}

\begin{proof}  The line bundle $\L = L(\Sigma)^2 \otimes K^{-1}$ is
$\widehat{LG^r}$-equivariant at levels $2(k+c),\ldots,2(k+c)$.  Since
$\M(\Sig)$ carries up to isomorphism a unique line bundle at every
level \cite[3.12]{me:co}, it follows that $\L$ is the pre-quantum line
bundle for the symplectic structure defined by $B_{2(k+c)}$.
Hence $\L_F$ is a pre-quantum line bundle for the corresponding
symplectic structure on $F$ (cf. \cite{al:fi}, Subsection 4.4.3), and $\int_F
\exp(\hh c_1(\L_F))$ is the symplectic volume $\Vol_{B_{k+c}}(F)$ for
the 2-form defined using $B_{k+c}$.
We claim that the symplectic volume coincides with the Riemannian
volume, which will complete the proof since 
$\Vol_{B_{k+c}}(T^2)=\# T_{k+c}$ 
by Lemma \ref{lem:bila} from Appendix \ref{app:groups}.  
By  our description of $M(\Sig_h^r)$ as a fusion product, 
the fixed point manifold $F=F(\Sig_h^r)$ is obtained from the
fixed point manifold $F(\Sig_0^2)$ (viewed as a group valued
Hamiltonian $T^2$-space) by fusion: $F(\Sig_1^1)=F(\Sig_0^2)_{\on{fus}}$
and
$$ F(\Sig_h^r)=F(\Sig_1^1)\fus\cdots \fus F(\Sig_1^1)\fus F(\Sig_0^2)
\fus\cdots \fus F(\Sig_0^2), $$
with $h$ factors $F(\Sig_1^1)$ and $(r-1)$ factors $F(\Sig_0^2)$.  
Lemma \ref{lem:fusabelian} from Appendix \ref{sec:fusion}
says that the symplectic volume of group valued Hamiltonian torus
spaces does not change under fusion. Hence 
$\Vol_{B_{k+c}}(F(\Sig_h^r))=\Vol_{B_{k+c}}(F(\Sig_0^2))^{h+r-1}$. 
Finally, the expression 
\eqref{eq:2form} for the 2-form on $M(\Sig_0^2)$ shows that 
$\Vol_{B_{k+c}}(F(\Sig_0^2))$ coincides with the 
Riemannian volume of $T^2$ with respect to $B_{k+c}$.  
\end{proof}

\begin{proposition} \label{prop:phase1} The phase factor is given by 
$\zeta_F(t_\lam)^{1/2} = (-1)^{(h+r-1)\,{\# \mf{R}_+}}.$
\end{proposition}

\begin{proof} 
The point $m = (e,\ldots,e) \in F$ lies in identity level set of
$\Phi$, and its stabilizer in $G^r$ is the image of the diagonal
embedding of $G$.  The 2-form $\om$ restricts to a symplectic form on
the tangent space $E=T_mM(\Sig)$.
%
%
By Equations (27) and (29) of \cite{al:fi}, $\zeta_F(t_\lam)^{1/2}$ can be
computed in terms of the symplectomorphism $A$ of $E$ defined by
$t_\lam$: Choose an $A$-invariant compatible complex structure on $E$
to view $A$ as a unitary transformation, and let $A^{1/2}$ be the
unique square root having all its eigenvalues in the set $\{e^{i\phi}
\, | \, 0\le \phi <\pi\}$. Then $\zeta_F(t_\lam)^{1/2}=\det(A^{1/2})$.

We first apply this recipe for the 2-holed sphere $\Sig_0^2$, 
so that $E=T_m M(\Sig_0^2)=\g\oplus\g$. 
Formula \eqref{eq:2form} shows that $\om_m$ is the standard
2-form on $\g\oplus\g$, given by the inner product $B$. A
compatible complex structure is given by the endomorphism 
$(\xi,\eta)\mapsto (-\eta,\xi)$. Thus, as a complex 
$G$-representation 
$E$ is just the complexification $E=\g^\C$. It follows 
that the eigenvalues of $A$ (other than $1$) come in complex 
conjugate pairs 
$$e^{i\phi_j},e^{-i\phi_j},\ \ \ 0<\phi_j\le\pi/2,$$
and the corresponding eigenvalues of $A^{1/2}$ are 
$e^{i\phi_j/2}$ and $e^{i\pi-i\phi_j/2} = - e^{-i\phi_j/2}$.
Hence
$$\zeta_F(t_\lam)^{1/2}=(-1)^{\#\mf{R}_+}.$$
Now consider the case $r \ge 1,h$ arbitrary.  The tangent space is
$T_mM(\Sig_h^r)=(\g\oplus\g)^{h+r-1}$, but because of the fusion terms
the symplectic form is not the standard symplectic form defined by the
inner product on $\g$.  However, by Appendix \ref{sec:fusion}, Lemma
\ref{lem:homotop} it is equivariantly and symplectically {\em
isotopic} to the standard symplectic form. Since the phase factor
$\zeta_F(t_\lam)^{1/2}$ is a root of unity, it is invariant
under equivariant symplectic isotopies, and we conclude as before that
$\zeta_F(t_\lam)^{1/2}=(-1)^{(h+r-1)\,{\# \mf{R}_+}}$.
\end{proof}

\subsection{Verlinde formula}
From Equation \eqref{eq:fp} we obtain, using 
\eqref{eq:d}, \eqref{eq:a} and 
Propositions \ref{prop:int} and \ref{prop:phase1}, 
\begin{theorem}[Verlinde Formula]\label{th:verlinde}
Let $G$ be a simply connected Lie group and $k$ a given integral 
level. The $\Spin_c$-index of the moduli space of flat connections 
on $\Sig_h^r$ at level $k$, with markings 
$\mu=(\mu_1,\ldots,\mu_r)\in(\Lambda^*_k)^r$ 
is given by the formula
\begin{equation}\label{eq:verlinde} 
\chi(\M(\Sig_h^r)_\mu)
={(\# T_{k+c})}^{h-1}
\sum_{\lambda\in\Lambda^*_k} |J(t_\lam)|^{2-2h}
\chi_{\mu_1}(t_\lam)^*\cdots  \chi_{\mu_r}(t_\lam)^*.
\end{equation}
\end{theorem}

\begin{remarks}
\begin{enumerate}
\item
Theorem \eqref{th:verlinde} also covers the case without boundary, since
$M(\Sig_h^0)=M(\Sig_h^1,0)$. One obtains
$$ \chi(\M(\Sig_h^0))={(\# T_{k+c})}^{h-1}
\sum_{\lambda\in\Lambda^*_k} |J(t_\lam)|^{2-2h}.
$$
\item 
For the two-holed sphere $\Sig_0^2$, formula \eqref{eq:verlinde}
simplifies by the orthogonality relations for level $k$ characters, 
and gives $\chi(M(\Sig_0^2)_{\mu_1,\mu_2})=\delta_{\mu_1,*\mu_2}$. 
\item
In Bismut-Labourie \cite{bi:sy} the $\Spin_c$-indices
$\chi(\M(\Sig_h^r)_\mu)$ are computed by 
direct application of the Kawasaki-Riemann-Roch formula to the 
reduced spaces. Their approach involves a description of all 
orbifold strata of the reduced space. 
The equality with the above sum over level $k$ weights is  
non-trivial; it is established in \cite{bi:sy} for 
sufficiently high level $k$.
\item
Theorem \ref{th:verlinde} gives a formula for a $\Spin_c$-index rather than
the dimension of a space of holomorphic sections. Vanishing results
for higher cohomology groups have recently been proved by Teleman
\cite[Section 8]{te:qu}.
\end{enumerate}
\end{remarks}
\vskip .1in
\section{Extension to non simply-connected groups}
\label{sec:prime}
In this section we consider moduli spaces of flat connections for
compact, connected semi-simple Lie groups that are not necessarily
simply connected, for surfaces with one boundary component.  The case
of multiple boundary components will be considered elsewhere.  Write
$G'=G/\Gamma$, where $G$ is simply connected and $\Gamma\subset Z(G)$
is a subgroup of the center $Z(G)$ of $G$. The covering $G\to G'$
identifies $\Gamma$ with the fundamental group $\pi_1(G')$.

The main problem in dealing with non-simply connected groups is that
the corresponding gauge groups are disconnected.  This happens already
for the loop group $LG'=\on{Map}(S^1,G')$: The kernel of the
natural map $LG'\to \pi_1(G')=\Gamma$ is the identity component
$L_0G'$ of $LG'$, identifying the group of components with $\Gamma$.

The disconnectedness of gauge groups causes a number of difficulties:
For example integrality of the level does not always guarantee the
existence of a pre-quantization of the moduli space.  In this section
we will obtain sufficient conditions for the existence of a
pre-quantization, and calculate the $\Spin_c$-index of moduli spaces
with prescribed holonomies around the boundaries for these cases.

\subsection{Gauge groups of surfaces for non-simply connected 
structure groups}
If $r\ge 1$, the surface $\Sig_h^r$ retracts onto a 
wedge of $2h+r-1$ circles. Hence every 
principal $G'$-bundle over $\Sig_h^r$ is trivial. 
On the other hand, principal $G'$-bundles over $\Sig_h^0$ 
are classified by elements of $\Gamma$: The bundle corresponding 
to $\gamma\in\Gamma$ is obtained by gluing the trivial 
bundles over $\Sig_h^1$ and over the disk $\Sig_0^1$, using 
a loop $g'\in LG'$ representing $\gamma$ as a transition 
function. 

Let $\Sig=\Sig_h^1$. As before we consider the gauge group
$\G'(\Sig)=\on{Map}(\Sig,G')$ and the subgroup $\G'_\p(\Sig)$ of gauge
transformation that are trivial on the boundary. To explain their
relationship with the gauge groups $\G(\Sig),\G_\p(\Sig)$ we consider
the homomorphism
\begin{equation}\label{eq:funda}
\G'(\Sig) \to \Gamma^{2h}
\end{equation}
which takes every gauge transformation $g':\,\Sig\to G'$ 
to the map induced on fundamental groups, viewed as an element 
of $\Hom(\pi_1(\Sig),\pi_1(G'))
=\Hom(\Z^{2h},\Gamma)=\Gamma^{2h}$. Let $\G(\Sig)\to\G'(\Sig)$ 
be the map induced by the covering $G\to G'$. 
\begin{proposition}\label{prop:exact}
The sequences
\begin{equation}\label{eq:based1}
1\to\Gamma\to \G(\Sig)\to \G'(\Sig)\to\Gamma^{2h}\to 1
\end{equation}
and 
\begin{equation} \label{eq:based}
1\to \G_\p(\Sig)\to\G'_\p(\Sig)\to\Gamma^{2h}\to 1
\end{equation}
are exact. Restriction to the boundary $\G'(\Sig)\to \G'(\p\Sig)\cong LG'$ 
takes values in $L_0G'$, giving rise to an exact sequence
\begin{equation} \label{eq:res}
1\to \G'_\p(\Sig)\to \G'(\Sig) \to L_0 G'\to 1.
\end{equation}
\end{proposition}
\begin{proof}
We first show that the map $\G'_\p(\Sig)\to \Gamma^{2h}$ is surjective.
Present $\Sig_h^0$ as a quotient of a $4h$-gon $P$, with 
sides identified according to the word 
$$C_1 C_2 C_1^{-1} C_2^{-1}\cdots C_{2h-1} C_{2h} C_{2h-1}^{-1} C_{2h}^{-1} .$$ 
Then $\Sig=\Sig_h^1$ is obtained as a similar quotient of $P$ minus a
disk in $\on{int}(P)$. The sides $C_j$ map to generators of
$\pi_1(\Sig)$, which we also denote $C_j$.  Given $\gamma\in
\Gamma^{2h}$, choose continuous maps $g_j':\,C_j\to G$ such that
$g_j'=e$ at the base point and such that the loop $g_j'$ represents
$\gamma_j$. Since $\pi_1(G')$ is abelian, the concatenation of loops
$$\prod_{j=1}^h [g'_{2j-1},g'_{2j}]:\,\prod_{j=1}^h
[C_{2j-1},C_{2j}]\to G'$$ 
is homotopically trivial. Hence the maps $g_j'$ extend to a continuous
map $g': \ \Sigma \to G'$ trivial on $\partial \Sigma$. By a
$C^0$-small perturbation, $g'$ can be changed to a smooth gauge
transformation, which still vanishes on $\p \Sig$.  We next show
exactness of \eqref{eq:based} at $\G'_\p(\Sig)$.  Since $G$ is simply
connected, a necessary and sufficient condition for an element
$g'\in\G'_\p(\Sig)$ to admit a lift $g\in \G_\p(\Sig)$ is that the
induced maps on fundamental groups be trivial. Thus $g'$ is in the
image of the map $\G_\p(\Sig)\to \G'_\p(\Sig)$ if and only if it is in
the kernel of the map $\G'_\p(\Sig)\to \Gamma^{2h}$.  Injectivity of
the map $\G_\p(\Sig)\to \G'_\p(\Sig)$ is obvious.  This proves
\eqref{eq:based} and the proof of \eqref{eq:based1} is similar.

Finally we prove \eqref{eq:res}.  Let $g'\in \G'(\Sig)$. Since
$\p\Sig$ is homotopic to $\prod_{j=1}^h [C_{2j-1},C_{2j}]$, the
restriction of $g'$ to the boundary defines a contractible loop in
$G'$.  Hence $g'|_{\p\Sig}\in L_0G'$.
\end{proof}

\subsection{Moduli spaces of flat connections}
View $\A(\Sig)=\Oma^1(\Sig,\g)$ as the space of $G'$-connections, 
and consider the action of the gauge group $\G'(\Sig)$. 
As before, the action of the subgroup $\G'_\p(\Sig)$ is Hamiltonian, with 
moment map the curvature.  The symplectic
quotient is the moduli space of flat $G'$-connections up to 
based gauge equivalence, 
$$ \M'(\Sig)=\A_{\on{{flat}}}(\Sig)/\G_\p'(\Sig).   
$$
It carries a residual Hamiltonian action of 
$L_0G'=\G'(\Sig)/\G'_\p(\Sig)$, with moment map induced by 
the pullback of a connection to the boundary. Using the 
surjection $LG\to L_0G'$, we will view $\M'(\Sig)$ as a 
Hamiltonian $LG$-space where $\Gamma\subset LG$ acts trivially.
We will also use the notation $\M(\Sig,G')=\M'(\Sig)$ to 
indicate the structure group. 

The moduli space $\M(\Sig)$ of flat $G$-connections is a finite
covering of $\M'(\Sig)$. Identify $\G_\p(\Sig)$ with the identity
component of $\G'_\p(\Sig)_0$. By Proposition \ref{prop:exact}, there
is an isomorphism
\begin{equation}
\G'(\Sig)/\G'_\p(\Sig)_0=\Gamma^{2h}\times \G'(\p\Sig)_0
\cong \Gamma^{2h}\times L_0G'. 
\end{equation}
This shows that the action of $LG$ on
$\M(\Sig)=\A_{\on{flat}}(\Sig)/\G_\p(\Sig)$ extends to an action of the
direct product $\Gamma^{2h}\times LG$, and
$$\M'(\Sig)=\M(\Sig)/\Gamma^{2h}.$$
Similarly, the holonomy manifold $M'(\Sig)$ of $\M'(\Sig)$ is just
$$ M'(\Sig)=M(\Sig)/\Gamma^{2h}=G^{2h}/\Gamma^{2h}=(G')^{2h}.$$
The product of commutators, $\Phi: \ M(\Sig)=G^{2h}\to G$ is invariant
under the action of $\Gamma^{2h}$ and descends to the $G$-valued
moment map $\Phi':\,M'(\Sig)\to G$.

For any $g' \in G'$, the moduli space of flat $G'$-connections with
holonomy conjugate to $g'$ is a disjoint union of symplectic quotients
$M'(\Sig)_g$ where $g$ varies over all pre-images of $g'$ in $G$.  The
reduced spaces at central elements $\gamma \in \Gamma\subset G$ may
also be interpreted as moduli spaces of flat connections on the
$G'$-bundle over $\Sig_h^0$, with topological type given by $\gamma$.

The moduli space $\M(\Sig,G')$ for the semi-simple group $G'$ 
is a finite cover of a product of moduli spaces for simple 
groups.  For $j=1,\ldots,s$, let
$\Gamma_j\subset G_j$ be the image of $\Gamma$ under projection to the
$j$th simple factor.  Then $G'$ covers the product of groups
$G_j'=G_j/\Gamma_j$, and since $\M(\Sig,\prod_{j=1}^s
G_j')=\prod_{j=1}^s \M(\Sig,G_j')$, one obtains a finite covering,
\begin{equation}\label{eq:coverings} 
\M(\Sig,G') \to \M(\Sig,G_1')\times \cdots\times \M(\Sig,G_s').
\end{equation}

\subsection{Fixed point manifolds}
\label{sec:fixed}
Since the $G$-action on $M'(\Sig)=(G')^{2h}$ is the conjugation 
action on each factor, the fixed point set for elements $g\in G$ is
the direct product of centralizers $G'_{g'}$ where $g'$ is the image 
of $g$. In contrast to the simply connected case, these 
centralizers may now be disconnected. For example, if $G'=\SO(3)$,
the centralizer of the group element given by rotation by $\pi$ about 
the $z$-axis is $\on{O}(2)\subset \SO(3)$. 
\begin{lemma}\label{lem:central} 
Let $t'$ be the image of a regular element $t\in T^{\on{reg}}$, and
$G'_{t'}$ its centralizer in $G'$.  Then $T'$ is a normal subgroup of
$G'_{t'}$ and $G'_{t'}/T'=W_{t'}$, the stabilizer group of $t'$ under
the action of the Weyl group.
\end{lemma}
\begin{proof} 
The image $g'\in G'$ of an element $g\in G$ is fixed under
$t'$ if and only if $[t,g]\in\Gamma$.  In this case, 
$\Ad_t$ fixes the maximal torus $\Ad_{g}(T)$ pointwise, 
hence  $\Ad_g(T)=T$ since $t$ is a regular element. 
Therefore $g\in N_G(T)$ and $g'\in N_{G'}(T')$.  
\end{proof}
The stabilizer group $W_{t'}$ can be re-interpreted as 
follows. Define an injective group homomorphism 
\begin{equation}\label{eq:imbed}
\varphi:\,Z(G)\to W,\, \gamma\mapsto w_\gamma
\end{equation} 
where $w_\gamma\in W$ is the unique element with
$\gamma\exp(\Alc)=w_\gamma(\exp(\Alc))$. 
Since the action of $Z(G)$ on $T$ commutes 
with the Weyl group action, it induces an action on the 
alcove $\Alc=T/W$.  The following is well-known:

\begin{lemma}\label{lem:isom}
For any $\xi\in\on{int}(\Alc)$, the homomorphism $\varphi$ restricts
to an isomorphism from the stabilizer group $\Gamma_\xi$ onto
$W_{t'}$, where $t'=\exp_{G'}(\xi)$. In particular, $W_{t'}$ is
abelian.
\end{lemma}

\begin{proof}
We describe the inverse map. Let $t=\exp_G(\xi)$.  Given $w\in
W_{t'}$, the element $\gamma=t^{-1}\, w(t)$ lies in $\Gamma$. The
equation $\gamma t=w(t)$ means, by definition of the action of
$\Gamma$ on $\Alc$, that $\gamma\in \Gamma_\xi$.
\end{proof}

From Lemmas \ref{lem:central} and \ref{lem:isom} we obtain the
following description of the fixed point sets.  For any $\gamma\in
\Gamma^{2h}$, let $F_\gamma\subset M'(\Sig)=(G')^{2h}$ be the
pre-image of $\varphi(\gamma)\in W^{2h}$ under the homomorphism
$$ (N_{G'}(T'))^{2h}\to (N_{G'}(T')/T')^{2h}=W^{2h}.$$
Then the fixed point set $M'(\Sig)^t$ for $t=\exp(\xi)$,
$\xi\in\on{int}(\Alc)$ is the union,
\begin{equation} \label{eq:union}
M'(\Sig)^t=\coprod_{\gamma\in\Gamma_\xi} F_\gamma.
\end{equation} 
Since the element $B_c^\sharp(\rho)\in\on{int}(\Alc)$ is fixed under
the action of $\Gamma$ (see \cite[Theorem 1.22]{bi:sy}), all of the
sub-manifolds $F_\gamma$ arise as fixed point manifolds of
$t=\exp(B_c^\sharp(\rho))$. In particular, $\om$ pulls back to a
symplectic form on each $F_\gamma$.
\begin{proposition}\label{prop:symplvol} 
For all $\gamma\in\Gamma^{2h}$, the symplectic volume of $F_\gamma$ is
equal to the Riemannian volume of $(T')^{2h}$ with respect to the
inner product $B$ defining $\om$.
\end{proposition}

\begin{proof}
We have $M'(\Sig_h^1)=M'(\Sig_1^1)\fus\ldots\fus M'(\Sig_1^1)$, and 
the fixed point set is just a fusion product of the fixed point sets of
the factors. Because of Lemma \ref{lem:fusabelian} 
it is enough to consider the case $h=1$. Thus
$M'(\Sig_1^1)$ is the $G$-valued Hamiltonian 
$G$-space $G'\times G'$ with
$G$ acting by conjugation, moment map $(a',b')\mapsto aba^{-1}b^{-1}$,
and 2-form 
\begin{equation} \label{eq:omdef}
\om_{(a',b')} =\hh\big( B(a^*\theta,b^*\olt)+B(a^*\olt,b^*\theta)
-B((ab)^*\theta, (a^{-1} b^{-1})^*\olt)\big)
\end{equation}
Using left-trivialization of the tangent bundle to identify
$T_{(a',b')}M'(\Sig_1^1) \cong\g\oplus\g$, and using the metric $B$ to
identify skew-symmetric 2-forms with skew-symmetric matrices, the
2-form $\om$ is given at $(a',b')$ by the block matrix,
\begin{equation} \label{eq:cdef}C=\left(\begin{array}{cc} C_{11}&C_{12}\\ C_{21}&C_{22}
\end{array}\right)\end{equation}
where $C_{ij}$ are the following endomorphisms of $\g$, 
\beq
C_{11}&=&\hh(\Ad_{b}-\Ad_{b^{-1}}),\\
C_{12}&=&\hh(-1+\Ad_{b}+\Ad_{a^{-1}}+\Ad_{ba^{-1}}),\\
C_{21}&=&\hh(1-\Ad_{a}-\Ad_{b^{-1}} -\Ad_{ab^{-1}}),\\
C_{22}&=&\hh(\Ad_{a^{-1}}-\Ad_{a} 
).  \eeq 

Suppose now that $F\subset G^2$ is the fixed point manifold labeled by
$(\gamma_1,\gamma_2)\in \Gamma^2$.  The 2-form $\om_F$ on $F$ is given
at any point $(a',b')\in F$ by the endomorphism $C^T=C|_{t\times\t}$,
with $\Ad_a,\Ad_b$ becoming the Weyl group action of
$w_1=\varphi(\gamma_1),\, w_2=\varphi(\gamma_2)$ on $\t$. Using that
$w_1,w_2$ commute one verifies that
$$ C^T_{11}C^T_{22}-C^T_{12}C^T_{21}=1,$$
showing $\det(C^T)=1$. Since the top degree part of $\exp(\om_F)$ is
equal to the standard volume form on $(T')^2$, times the Pfaffian of
$C^T$, this shows $\Vol(F)=\Vol(T')^2$.
\end{proof}

Later we will need the following remarkable fact.

\begin{lemma}\label{lem:meet}
The identity level set $(\Phi')^{-1}(e)$ intersects each of the
sub-manifolds $F_\gamma$.
\end{lemma}

\begin{proof}
It suffices to consider the case of the 1-punctured torus
$\Sig=\Sig_1^1$. Let $\gamma=(\gamma_1,\gamma_2)$ be given.  By
Proposition \ref{prop:lifts} from Appendix \ref{app:lifts}, it is
possible to choose commuting lifts $g_1,g_2\in N_G(T)$ of
$w_1=\varphi(\gamma_1),\ w_2=\varphi(\gamma_2)$. Then $(g_1',g_2')\in
F_\gamma$ and $\Phi'(g_1',g_2')=[g_1,g_2]=e$.
\end{proof}

\section{Pre-quantization of $\M'(\Sig)$}
\label{sec:preprime}

The purpose of this section is to prove the following theorem.
\begin{theorem}\label{th:classline}
\begin{enumerate}
\item If $L' \to \M'(\Sig)$ is an $\LG$-equivariant line bundle at
level $l$, then $L'$ admits a connection whose curvature is the
two-form defined by $B_l$.
\item Any two $\LG$-equivariant level $l$ line bundles differ by the
pull-back of a flat line bundle on $(G')^{2h}$ defined by a character
of $\Gamma^{2h}$.
\item $\M'(\Sig)$ admits a level $l=(l_1,\ldots,l_s)$ line bundle if
each $l_j$ is a multiple of the greatest common divisor of $2c_j$ and
$\#\Gamma_j^2$.
\item For each such $l$, there is a unique level $l$ line bundle with
the property that for all $t\in T^{\on{reg}}$, and all fixed points
$\hat{m}'\in \M'(\Sig)^t$ in the zero level set of the moment map, the
action of $t$ on the fiber at $\hat{m}'$ is trivial.
\end{enumerate}
\end{theorem}

\begin{example}
We consider the adjoint groups $G'=G/Z(G)$ for the simple
groups. From the table in Appendix \ref{app:groups}
we obtain the following smallest positive levels $l_0$
for which our criterion guarantees a level $l_0$ line bundle 
over $\M'(\Sig)$. 

\begin{tabular}{c@{ : }l}
$A_N'$ & $l_0=N+1$ for $N$ even, and $l_0=2(N+1)$ for $N$ odd. \\
$B_N'$ & $l_0=2$.\\
$C_N'$ &$l_0=2$ for $N$ even, $l_0=4$ for $N$ odd.\\
$D_N'$ &$l_0=4$ for $N$ even, $l_0=16$ for $N=5,9,\ldots$,
$l_0=8$ for $N=7,11,\ldots$.\\
$E_6'$& $l_0=3$.\\
$E_7'$& $l_0=4$.
\end{tabular}

\end{example}

The condition for the adjoint group of $\SU(n)=A_{n+1}$ coincides with
the condition from Beauville-Laszlo-Sorger \cite{be:ve,be:pg}. 
Recent results of Fuchs--Schweigert \cite{fu:ac} (also Teleman
\cite{te:qu}) suggest that these conditions are not optimal.

\begin{proof}[Proof of Theorem \ref{th:classline}]
Recall that $\M(\Sig)$ carries a unique equivariant line bundle
$L^{(l)}(\Sig)$ at level $l$ \cite[Proposition 3.12]{me:co}.
Moreover, this line bundle carries a connection $\nabla^{(l)}$ whose
curvature is the two-form defined by $B_l$.  The pull-back of $L'$ to
$\M(\Sig)$ is isomorphic to $L^{(l)}(\Sig)$.  The average of
$\nabla^{(l)}$ over $\Gamma^{2h}$ descends to a connection on
$L'$.

For part (b), we may assume $L'$ is a line bundle at level $l=0$.  The
quotient $L'/\Omega G$ is a $G$-equivariant line bundle over
$M'(\Sig)=(G')^{2h}$. Since $H^2((G')^{2h},\R)=
H^2_G((G')^{2h},\R)=0$, every $G$-line bundle over $M'(\Sig)$ is
isomorphic to a flat line bundle defined by a character of
$\Gamma^{2h}$ (see e.g. \cite{mu:li}).

The discussion in Subsection \ref{sec:canon} (resp.
\ref{sec:central}) shows that there exists a level $l_j$ line bundle over
$\M(\Sig,G_j')$ if $l_j$ is a multiple of $2c_j$ (resp.
$\#\Gamma_j^2$) with the property given in part (d) of the theorem.
The theorem follows by taking tensor products and pull-backs under the map
\eqref{eq:coverings}.
\end{proof}

\subsection{Canonical bundles}
\label{sec:canon}

The anti-canonical line bundles over each factor $\M(\Sig,G_j')$ in
the decomposition \eqref{eq:coverings} are at level $2c_j$.  It
remains to show that for any $t \in T^{\on{reg}}$ and all fixed points
$\hat{m}'\in \M'(\Sig,G_j)^t$ in the zero level set of the moment map
, the eigenvalue $\kappa(\hat{m}',t)$ for the action of $t$ on the
fiber at $\hat{m}'$ is trivial.  We will prove a  
stronger
statement, which we will need later in the fixed point theorem.
Recall from \cite{al:fi}, Section 4.4 the definition of the square root
$\kappa(\hat{m}',t)^{1/2}$.

\begin{proposition} \label{prop:sqrt} 
The square root for the eigenvalue of $t$ on the fiber of the 
canonical line bundle at $\hat{m}'$ is given by 
$\kappa(\hat{m}',t)^{1/2} = (-1)^{h\# \mf{R}_+} .$
\end{proposition}

\begin{proof}  
Let $m'\in (\Phi')^{-1}(e)$ be the image of $\hat{m}'$ in
$M'(\Sig)^t$, and $F_\gamma$ the connected component containing
$m'$. 
By (29) of \cite{al:fi}, the square root $\kappa(\hat{m}',t)^{1/2}$ coincides
with the square root of the eigenvalue of the action of $t$ on the
symplectic vector space $T_{m'}M'(\Sig)$.

As in the simply connected case, it is enough to consider the case of
the 1-punctured torus $\Sig=\Sig_1^1$. We adopt the notation and
definitions \eqref{eq:omdef},\eqref{eq:cdef} from the proof of Proposition
\ref{prop:symplvol}.

Suppose $\k\subset \g$ is invariant under both of the commuting
operators $\Ad_a$ and $\Ad_b$. From the formula for $\om_{(a',b')}$ it
follows that
$$\g\oplus\g=(\k\oplus\k)\oplus(\k^\perp\oplus \k^\perp)$$
as an $\om$-orthogonal direct sum.  Take $\k$ to be the direct sum of
joint eigenspaces for $\Ad_a,\Ad_b$ with eigenvalues $(-1,1)$,
$(-1,-1)$ or $(1,-1)$.  Then  
$$C_{11} |_\k =C_{22} |_\k=0, \ \ \ 
C_{12} |_\k =-\on{Id}=-C_{21} |_\k .$$ 
Hence $\om$ restricts to {\em minus} the standard symplectic form on
$\k\oplus\k$.  Therefore the square root for the action on the
canonical bundle ${\det}_\C^{-1}(\k\oplus\k)$ is
$(-1)^{(\dim\k-\dim(\k)^t)/2}$.

Consider on the other hand the $\k^\perp\subset\g$ contribution.  We
claim that the symplectic form $C$ on $\k^\perp\oplus\k^\perp$ is
homotopic to the standard symplectic form $C_0=\left(\begin{array}{cc}
0&\on{Id}\\-\on{Id}&0\end{array}\right)$, by a homotopy of symplectic
forms which are invariant under the action of $t$.  This will
imply that the square root for the action on
${\det}_\C^{-1}(\k^\perp\oplus\k^\perp)$ is $
(-1)^{(\dim\k^\perp-\dim(\k^{\perp})^t)/2}$, and finally
$$ \kappa(\hat{m}',t)^{1/2} = (-1)^{(\dim\g-\dim\t)/2}=(-1)^{\# \mf{R}_+} .$$
For the rest of this paragraph we consider the restriction $C^\perp$
of $C$ to $\k^\perp\oplus\k^\perp$.  For $s\in[0,1]$ let
$C^\perp_s=(1-s)C^\perp_0+s C^\perp$.  It suffices to show that $\det(C^\perp_s)\not=0$ for
$s\in [0,1]$.  We calculate
\beq 
\det(C^\perp_s)&=&
\det\left(\begin{array}{cc}
s C^\perp_{11}& 1+s(C^\perp_{12}-1)\\-1+s(C^\perp_{21}+1)&s C^\perp_{22}
\end{array}\right)
 \\
&=&\det\big(s^2(C^\perp_{11}C^\perp_{22}-C^\perp_{12}C^\perp_{21}
-C^\perp_{12}+C^\perp_{21}+1
)-s(C^\perp_{21}-C^\perp_{12}+2)
+1\big)\\
&=&\det\big((s^2-s)(2+(C^\perp_{21}-C^\perp_{12}))+1 \big)
\eeq
For $s\in[0,1]$ one has $0\ge s^2-s\ge -\f{1}{4}
$. Hence it suffices to show that all eigenvalues of 
the symmetric operator $C^\perp_{21}-C^\perp_{12}$ are strictly 
smaller than $2$. But 
$$ C^\perp_{21}-C^\perp_{12}=1-\hh(\Ad_{a}+\Ad_{a^{-1}}+
\Ad_{b}+\Ad_{b^{-1}}+\Ad_{a\,b^{-1}}+     
\Ad_{b\,a^{-1}}).
$$
%
On the joint eigenspace of $\Ad_a,\Ad_b$ with eigenvalue $(1,1)$, this
is strictly negative.  On the orthogonal complement of the eigenspaces
for eigenvalue pairs $(-1,-1), (-1,1), (1,-1)$ and $(1,1)$, the operators
$\Ad_{a},\Ad_{b}$ are represented by 2-dimensional rotations by angles
$\psi_a,\psi_b$ not both of which are multiples of $\pi$. On any such
2-plane, the operator $C^\perp_{21}-C^\perp_{12}$ becomes
$(1-\cos(\psi_a)-\cos(\psi_b)-\cos(\psi_a-\psi_b))\on{Id}.$ The claim
follows since $1-\cos(\psi_a)-\cos(\psi_b)-\cos(\psi_a-\psi_b)<2$.
\end{proof}

\subsection{Central extension of the gauge group $\G'(\Sig)$}
\label{sec:central}

In order to define a pre-quantum bundle over $\M(\Sig,G')$ at level
$\# \Gamma^2$, we imitate the construction from the simply-connected
case.  The definition \eqref{eq:cocycle} of the cocycle carries over
to the non-simply connected case, and defines a central extension
$\wh{\G}'(\Sig)$ of the gauge group. 
\begin{proposition} \label{prop:divis}
Suppose each $k_j$ is a multiple of $\# \Gamma^2$.  Then there is a
canonical trivialization of the central extension $\wh{\G}'(\Sig)$
over $\G'_\p(\Sig)$.
\end{proposition}
\begin{proof}
Let $\ol{\Sig}$ be the surface obtained from $\Sig$ by capping off the
boundary component.  Let $\G'_c(\Sig)\subset \G'_\p(\Sig)$ 
denote the kernel of the
restriction map $ \G'(\ol{\Sig}) \to \G'(\ol{\Sig} -\Sig)$.  As
explained in \cite[p. 431]{me:lo}, it suffices to show that the
cocycle is trivial over the subgroup $\G'_c(\Sig)$ of $\G'_\p(\Sig)$.
We would like to define a map $ \alpha:\,\G_c'(\Sig)\to \U(1)$ with
coboundary condition
\begin{equation}\label{eq:coboundary}
\alpha(g_1'g_2')=\alpha(g_1')\alpha(g_2')c(g_1',g_2')
\end{equation}

Given $g'\in \G'_c(\Sig)$, let 
$\gamma\in\Gamma^{2h}\cong \Hom(\pi_1(\Sig),\Gamma)$ be its image
under the map \eqref{eq:funda}. The element $\gamma$  defines a 
covering $\Sig^\gamma$ of $\Sig$ by $\Gamma$, which is a possibly 
disconnected surface with $\# \Gamma$ boundary components. 
Choose a base point of
$\Sigma^\gamma$ mapping to the base point of $\Sigma$, and let
$\pi^\gamma:\,\Sig^\gamma\to \Sig$ be the covering projection.  The
pull-back $(\pi^\gamma)^*g'$ admits a unique lift $g\in
\G(\Sig^\gamma)$, with $g=e$ at the base point and $g(\lambda\cdot
x)=\lambda g(x)$ for all $\lambda\in\Gamma$, $x\in \Sig^\gamma$.
Note that $g$ is constant along each of the boundary circles of 
$\Sig^\gamma$.

The covering $\Sig^\gamma\to \Sig$ extends to a covering
$\ol{\Sig}^{\gamma}\to\ol{\Sig}$ over the capped-off surface.  Extend
$g \in \G(\Sig^\gamma)$ first to $\ol{\Sig}^{\gamma}$ by the constant
extension on the capping disks, and then further to a map $\ol{g} \in
\G(\Sig^\gamma\times [0,1])$, in such a way that the extension is
trivial on $\ol{\Sig}^\gamma\times \{1\}$.  Define
\begin{equation} \label{alpha}
\alpha(g')=\exp(2\pi i \f{1}{{\# \Gamma}} \int_{\ol{\Sig}^\gamma \times [0,1]}
\ol{g}^*\eta).
\end{equation}
Since 
\begin{equation} \label{eq:whatweneed}
\f{1}{\#\Gamma^2} [\eta]\in  H^2(G,\Z)
\end{equation}
by our assumption on the level, the definition \eqref{alpha} is
independent of the choice of $\ol{g}$.  It remains to verify the
coboundary condition.  Let $g_1',g_2'\in\G_c'(\Sig)$ and
$\gamma_1,\gamma_2\in \Gamma^{2h}$ their images. Define
$\ol{\Sigma}^{\gamma_1,\gamma_2}$ as the fibered product of
$\ol{\Sigma}^{\gamma_1}$ and $\ol{\Sigma}^{\gamma_2}$.  Let
$\pi_1,\pi_2,\pi_{12}$ denote the projections onto
$\ol{\Sigma}^{\gamma_1},\ol{\Sigma}^{\gamma_2}, \ol{\Sigma}^{\gamma_1
\gamma_2}$ respectively.  Let $g_{12}' = g_1' g_2' \in \G_c'(\Sigma)$.
We have canonical lifts $g_1\in \G(\Sigma^{\gamma_1}), g_2\in
\G(\Sigma^{\gamma_2})$ and $g_{12} \in\G(\Sigma^{\gamma_1 \gamma_2})$.
Choose extensions $\ol{g}_1,\ol{g}_2$ and $\ol{g_{12}}$ as above. Then
both $\pi_{12}^*\ol{g_{12}}$ and $ \pi_{1}^* \ol{g_1}\ \pi_{2}^*
\ol{g_2}$ are equal on $\ol{\Sig}^{\gamma_1,\gamma_2} \times \{ 0 \}$. Using
\eqref{eq:whatweneed},
\beq
\frac{1}{\# \Gamma} \int_{\ol{\Sigma}^{\gamma_1\gamma_2} \times [0,1]}
\ol{g_{12}}^* \eta &=&
\frac{1}{\# \Gamma^2}
\int_{\Sigma^{\gamma_1,\gamma_2} \times [0,1]} 
(\pi_{12}^*\ol{g_{12}})^* \eta \\
&=& \frac{1}{\# \Gamma^2}
\int_{\Sigma^{\gamma_1,\gamma_2} \times [0,1]} ( \pi_{1}^* \ol{g_1}\ 
\pi_{2}^* \ol{g_2})^* \eta \ \ \ \on{mod}\Z
\eeq
For the last term we compute, using the property of the 3-form $\eta$
under group multiplication $\Mult_G:\,G\times G\to G$, 
$$\Mult_G^*\eta=\pr_1^*\eta+\pr_2^*\eta+
\d \big(\hh B(\pr_1^*\theta,\pr_2^*\olt)\big)$$
(where $\pr_j:\,G\times G\to G$ are projections 
to the respective $G$-factor),
\begin{eqnarray*}
\lefteqn{
\frac{1}{\# \Gamma^2}
\int_{\ol{\Sigma}^{\gamma_1,\gamma_2} \times [0,1]} ( \pi_{1}^* \ol{g_1}\ 
\pi_{2}^* \ol{g_2})^* \eta \ \ \ 
}\\  
&=& \frac{1}{\# \Gamma^2}
\int_{\ol{\Sigma}^{\gamma_1,\gamma_2} \times [0,1]} 
\big(\pi_{1}^* \ol{g_1}^* \eta + \pi_{2}^* \ol{g_2}^* \eta + 
\hh \d B( \pi_{1}^*\ol{g_1}^*\theta , \pi_{2}^* \ol{g_2}^*
\ol{\theta} ) \big)\\ 
&=& \frac{1}{\# \Gamma} \int_{\ol{\Sigma}^{\gamma_1}
\times [0,1]} \pi_{1}^* \ol{g_1}^* \eta + \frac{1}{\#
\Gamma} \int_{\ol{\Sigma}^{\gamma_2} \times [0,1]} \pi_{2}^* \ol{g_2}^* \eta 
 + \hh \int_{\Sigma} B( (g_1')^*
\theta , (g_2')^* \ol{\theta}).
\end{eqnarray*}
This shows that the cocycle condition holds and completes the proof.
\end{proof}
As a consequence, if $k$ is a multiple of $\#\Gamma^2$ the moduli 
space $\M'(\Sig)$ carries a line bundle 
$$L'(\Sig)=(\A_{\on{flat}}(\Sig)\times\C)/\G'_\p(\Sig)$$ 
equipped with an action of $\wh{\G}'(\Sig)/\G'_\p(\Sig)$ where the
central circle acts with weight 1. The map $LG\to L_0G'
={\G}'(\Sig)/\G'_\p(\Sig)$ induces a map from the central extension at
level $k$, $\wh{LG}\to \wh{\G}'(\Sig)/\G'_\p(\Sig)$ which restricts to
the identity map on the central circle and has kernel $\Gamma$.  The
upshot is that $L'(\Sig)$ is a level $k$ line bundle over $\M'(\Sig)$,
where $\Gamma\subset G\subset \wh{LG}$ acts trivially.  It remains to
verify property (d) of Theorem \ref{th:classline}.
\begin{lemma}\label{lem:trivia}
Let $t\in T^{\on{reg}}\subset LG$ a regular element, and
$\hat{m}'\in\M'(\Sig)$ a fixed point for the action of $t$ in the zero
level set of the moment map. Then $t$ acts trivially on the fiber
$L'(\Sig)_{\hat{m}'}$.
\end{lemma}

\begin{proof}  If $t$ is contained in the identity 
component of the stabilizer $LG_{\hat{m}'}$, then the formula follows
from the pre-quantum condition.

Our strategy is to reduce to this case, using finite covers.  Since
$t$ fixes $\hat{m}'$, there exists a flat connection $A$ mapping to
$\hat{m}'$ and a gauge transformation $g' \in\G'(\Sig)$ restricting to
$t'$ on the boundary such that $g'\cdot A=A$.  The eigenvalue for the
action of $t$ on $L'(\Sig)$ is equal to the eigenvalue for the action
of $g'$ on the fiber $\{A\}\times \C$ of the level $k$ pre-quantum
line bundle over $\A(\Sig)$.

Let $\gamma\in \Gamma^{2h}$ be the image of $g'$ under the map
$\G'(\Sig)\to \Gamma^{2h}$.  As in the proof of Proposition
\ref{prop:divis} let $\pi^\gamma:\,\Sig^\gamma\to \Sig$ the covering
defined by $\gamma$, and $g\in \G(\Sig^\gamma)$ be the lift of the
pull-back $(\pi^\gamma)^*g'$. Clearly $g$ fixes
$A^\gamma:=(\pi^\gamma)^*A$.

The pull-back map $\G'(\Sig)\to \G'(\Sig^\gamma)$ lifts to the 
central extensions if one changes the level. Indeed, the 
cocycle for the central extension at level $k$ pulls back to the 
cocycle for the extension at level $k/\#\Gamma$. The map is 
compatible with the given trivializations over $\G_\p'(\Sig)$ and 
$\G_\p'(\Sig^\gamma)$, hence they define a commutative diagram
\begin{equation}
\vcenter{\xymatrix{
\G_\p'(\Sig)\ar[d]  \ar[r]& \G_\p'(\Sig^\gamma) \ar[d]  \\
         \wh{\G}'(\Sig)^{(k)}     \ar[r] & \wh{\G}'(\Sig^\gamma)^{(\f{k}
{\#\Gamma})}
}}
\end{equation}
It follows that the eigenvalue of $g'$ on the fiber $\{A\}\times\C$ is
equal to that of $g$ on the fiber $\{A^\gamma\}\times\C$.  By Lemma
\ref{lem:identitycomp} below, $g$ is contained in the identity
component of $\G(\Sig^\gamma)$. Hence its action on the line-bundle is
determined by the moment map (Kostant's formula), and hence is trivial
since $A^\gamma$ pulls back to the zero connection on $\p\Sig$.
\end{proof}

In the proof we used the  following Lemma.
\begin{lemma}\label{lem:identitycomp}
Let $\Sig_h^r$ be any oriented surface (possibly with boundary), and
$g\in\G(\Sig_h^r)$, $A\in\A_{\on{flat}}(\Sig_h^r)$ with $g\cdot
A=A$. Suppose that for some point $x\in \Sig_h^r$, $g(x)$ is a regular
element.  Then $g$ is contained in the identity component of the
stabilizer $\G(\Sig_h^r)_A$.
\end{lemma}

\begin{proof}
It is well-known that evaluation at $x$, $\G(\Sig_h^r)\to G$ restricts
to an injective map $\G(\Sig_h^r)_A\hra G$. The image $G_A$ of this
map is the centralizer of the image $\Hol_A$ of the homomorphism
$\pi_1(\Sig_h^r,x)\to G$ defined by parallel transport using $A$.
Since $ \Hol_A \subset G_{g(x)}$ and $G_{g(x)}$ is a torus, $G_{g(x)}
\subset G_A$.  It follows that $g$ is contained in the identity
component of $\G(\Sig_h^r)_A$.
\end{proof}

\subsection{The phase factor}

Using Theorem \ref{th:classline} and Proposition \ref{prop:sqrt}
we compute the phase factor appearing in the fixed point formula.
\begin{proposition} \label{prop:phase} 
Let $k = (k_1,\ldots,k_s)$, where each $k_j$ is a positive multiple of
the greatest common divisor of $c_j$ and $\#\Gamma_j^2$. Let $L'(\Sig)
\to \M'(\Sig)$ be the pre-quantum line bundle at level $k$ corresponding
to a character $\phi: \ \Gamma^{2h} \to U(1)$
(cf. Theorem \ref{th:classline}).  Let $t\in T^{\on{reg}}$,
and $F_\gamma \subset M'(\Sigma)$ any fixed point component.  The
phase factor $\zeta_{F_\gamma}(t)^{1/2}$ for the action on the
$\L_{F_\gamma}$ is given by
$$\zeta_{F_\gamma}(t)^{1/2}= \phi(\gamma) (-1)^{h\#\mf{R}_+}.$$
\end{proposition}

\begin{proof}
Let $m'\in F_\gamma\cap \Phi^{-1}(e)$, and $\hat{m}'\in \Phi^{-1}(0)$
the unique element projecting to $m'$. By Theorem \ref{th:classline},
the weight for the action of $t$ on $L'(\Sig)_{\hat{m}'}$ is
$\phi(\gamma)$.  On the other hand, by Proposition \ref{prop:sqrt},
$\kappa(\hat{m}',t)^{1/2} = (-1)^{h\# \mf{R}_+}$, which completes the
proof.
\end{proof}

\section{Verlinde formulas for non simply-connected groups}
\label{sec:verprime}

In this section we apply the fixed point formula to the Hamiltonian
$LG$-manifold $\M'(\Sig)=\M(\Sig,G')$.  The main result of the paper is
\begin{theorem}
 \label{th:ns}
Let $k = (k_1,\ldots,k_s)$, where each $k_j$ is a positive multiple of
the greatest common divisor of $c_j$ and $\#\Gamma_j^2$.  Let
$L'(\Sig)\to\M'(\Sig)$ be the pre-quantum line bundle at level $k$
corresponding to a character $\phi\in \Hom(\Gamma^{2h},\U(1))$, as in
Theorem \ref{th:classline}.  For any $\mu\in\Lambda^*_k$, the
$\on{Spin}_c$-index of the symplectic quotient $\M'(\Sig)_\mu$ at
$\mu$ is given by the formula,
\begin{equation} \label{eq:nsc} 
\chi(\M'(\Sig)_\mu) = 
\f{(\# T_{k+c})^{h-1}}{\# \Gamma^{2h}}
\sum_{ \lam \in \Lambda^*_k} \eps(\phi,\lambda)\
\frac{\# \Gamma_{\lambda}^{2h}}
{|J(t_\lambda)|^{2h-2}}\ \chi_{\mu}(t_\lambda)^*
.\end{equation}
Here $\eps(\phi,\lambda)=1$ if $\phi$ restricts to the trivial
homomorphism on $\Gamma_{\lambda}^{2h}$, and $0$ otherwise.
\end{theorem}

\begin{proof} 
Recall $t_\lambda = \exp(B_{k+c}^\sharp(\lambda+\rho))$.  Since
$B_c^\sharp(\rho)$ is fixed by $\Gamma$, the stabilizer of
$B_{k+c}^\sharp(\lambda+\rho)$ is equal to the stabilizer of
$\lambda$, using the embedding $\Lambda^*_k \to \Alc$ induced by
$B_k$.  By \eqref{eq:union}, the fixed point components for the action
of $t_\lambda$ on $M'(\Sig)$ are the sub-manifolds $F_\gamma$ with
$\gamma\in\Gamma_{\lambda}^{2h}$.  Since $F_\gamma$ is diffeomorphic
to $T^{2h}$, we have $\hat{A}(F_\gamma)=1$.  The normal bundle
$\nu_{F_\gamma}$ is $t_\lambda$-equivariantly isomorphic to the
constant bundle with fiber $(\g/\t)^{2h}$.  Indeed, left translation
by any $m=(g_1',\ldots,g_{2h}')\in F_\gamma$ on $(G')^{2h}$ commutes
with the action of $t_\lam$, and induces a $t_\lam$-equivariant
isomorphism of $\nu_{F_\gamma}(m)$ with $(\g/\t)^{2h}$. Therefore,
\begin{equation}\label{eq:DR}
\D_\R(\nu_{F_\gamma},t_\lam) = J(t_\lam)^{2h}=
(-1)^{h\#\mf{R}_+} \,|J(t_\lam)|^{2h}.
\end{equation}
Since the line bundle $\L$ is pre-quantum at level $2(k+c)$, the
integral $\int_{F_\gamma} \exp(\hh c_1(\L_F))$ is equal to the
symplectic volume with respect to the inner product $B_{k+c}$.  By
Proposition \ref{prop:symplvol} the symplectic volume is equal to the
Riemannian volume. Since $T'=T/\Gamma$, this shows
\begin{equation}
\int_{F_\gamma} 
\exp(\hh c_1(\L_{F_\gamma}))=
\f{\Vol_{B_{k+c}}(T^{2h})}{(\#\Gamma)^{2h}}
=\f{(\# T_{k+c})^{h}}{(\#\Gamma)^{2h}}.
\end{equation} 
Miraculously, the contribution is independent of $\gamma$.  The fixed point formula hence gives the following expression for 
the $\Spin_c$-index of a symplectic quotient $\M'(\Sig)_\mu$, 
for $\mu\in\Lambda^*_k$:
\begin{equation}\label{eq:completlygeneral}
\chi(\M'(\Sig)_\mu)=
\f{(\# T_{k+c})^{h-1}}{(\# \Gamma)^{2h}}
\sum_{\lambda\in \Lambda^*_k}
\f{\chi_\mu(t_\lam)^*}{|J(t_\lam)|^{2h-2}}\ 
(-1)^{h\#\mf{R}_+} \sum_{\gamma\in\Gamma_{\lambda}^{2h}}
\zeta_{F_\gamma}(t_\lam)^{1/2}.
\end{equation}
By the computation of the phase factor
in Proposition \ref{prop:phase}, 
$$\sum_{\gamma\in\Gamma_{\lambda}^{2h}}
\zeta_{F_\gamma}(t_\lam)^{1/2} 
=
(-1)^{h\#\mf{R}_+} \sum_{\gamma\in\Gamma_{\lambda}^{2h}} \phi(\gamma)
= (-1)^{h\#\mf{R}_+} \eps(\phi,\lambda)\
\# \Gamma_{\lambda}^{2h} $$
which completes the proof.
\end{proof}

It is remarkable that the right-hand side of \eqref{eq:nsc} always
produces  an 
integer, under the given assumptions on $k$.

\begin{remark}[Contribution from the exceptional element] 
\label{special element}
Suppose $G$ is simple.  As mentioned above, the element
$\xi_0=B_c^\sharp(\rho)$ is always a fixed point for the action of
$\Gamma$. The element $\lambda_0=B_k^\flat(\xi_0)=\f{k}{c}\rho$ is a
weight if and only if $c|k$. By \cite[Proposition 1.10]{bi:sy}, we
have $|J(\exp\xi_0)|^2=\# T_c=(1+\f{k}{c})^{-\dim T}\# T_{k+c}$.
Hence the contribution of $\lambda_0$ is
$$ (1+\f{k}{c})^{(h-1)\dim T} \eps(\phi,\lambda_0) \chi_\mu(t_{\lam_0})^*.
$$
Since $\Gamma_{\lambda_0}=\Gamma$, the factor $\eps(\phi,\lambda_0)$ is 
$1$ if $\phi$ is trivial and $0$ if $\phi$ is non-trivial. A theorem of 
Kostant \cite[Theorem 3.1]{ko:ma} asserts that the character value 
$\chi_\mu(t_{\lambda_0})$ is either $-1,0$, or $1$. More precisely, 
Kostant proves \cite[Lemma 3.6]{ko:ma} that
$$\chi_\mu(t_{\lambda_0})=\left\{\begin{array}{cl} 
(-1)^{l(w)}&\mbox{ if }\exists w\in W:\ w(\mu+\rho)-\rho\in \ol{R}'\\
0&\mbox{ otherwise }
\end{array}\right.
$$
where $\ol{R}'\subset\ol{R}$ is the lattice generated by all $c m_\alpha \alpha$, 
where $m_\alpha=||\alpha_0||^2/||\alpha||^2 \in \{1,2,3\}$.  
In particular, for simply laced groups $\ol{R}'=c\ol{R}$.
\end{remark}
\subsection{The case $G' = \on{PSU}(p)$.}
Let us specialize the above formulas to $G=\SU(p)$ and 
$\Gamma=Z(G)=\Z_p$ so that $G' = \on{PSU}(p)$. 
In this case Theorem \ref{th:ns} reduces   to the formulas of
Pantev \cite{pa:cm} ($p = 2$), and Beauville \cite{be:ve} ($p\ge 3$
prime), as follows. Recall $c=p$ and that for $p$ odd, $\M'(\Sig)$ is 
pre-quantizable for $k$ any multiple of $p$. The alcove 
$\Alc$ is a simplex with $p$ vertices, and $\Gamma=\Z_p$ acts by
rotation inducing a cyclic permutation of the vertices.  If $p$ is
prime, the only point of $\Alc$ with non-trivial 
stabilizer is the center $\xi_0=B_c^\sharp(\rho)$ of the alcove. 
By Remark \ref{special element}, the fixed point contribution for 
the exceptional element is 
$$ (1+\f{k}{p})^{(p-1)(h-1)}\eps(\phi,\lambda_0)\ \chi_\mu(t_{\lam_0})^*.$$  
Let $\chi(\M(\Sig)_\mu)$ be the index for the moduli space
of $G$-connections. Its fixed point contributions are exactly as
for $M'(\Sig)_\mu$, except for the overall factor $(\# \Gamma)^{-2h}
=p^{-2h}$ and a different weight for the exceptional element $t_{\lambda_0}$.  This shows,
$$ \chi(\M'(\Sig)_\mu)=p^{-2h}\chi(\M(\Sig)_\mu)
+(\eps(\phi,\lambda_0)-p^{-2h}) 
(1+\f{k}{p})^{(p-1)(h-1)}\chi_\mu(t_{\lambda_0})^*,$$
which (for $\mu=0$ and $\phi=1$) is exactly the formula given by Beauville.

\subsection{The sum over components}
Let $\mu\in\Alc$, and 
$\Co'_\mu \subset G'$  the conjugacy class of $\exp_{G'}(\mu)$.
The moduli space $\M'(\Sig,\Co'_\mu)$ of flat $G'$-connections 
on $\Sig=\Sig_h^1$, with 
holonomy around the boundary in $\Co'_\mu$ is a disjoint union
$$ \M'(\Sig,\Co'_\mu)=
\coprod_{\gamma \in \Gamma/\Gamma_\mu} \M'(\Sig)_{\gamma \mu}.$$
(For $\mu = 0$, the components may be identified with flat
$G'$-bundles over the closed surface $\Sig_h^0$.)  
The action of $\Gamma$ on $\Alc$ preserves the set $\Lambda^*_k\subset
\Alc$ of level $k$ weights. If $\mu \in \Lambda_k^*$, the
$\on{Spin}_c$-indices for the space $\M'((\Sig)_{\gamma \mu})$ can be
computed using Theorem \ref{th:ns}.  The sum over $\Gamma/\Gamma_\mu$
can be carried out using the following transformation property
of level $k$ characters
(see Bismut-Labourie \cite[Theorem 1.33]{bi:sy})
\begin{equation}\label{eq:trafo}
\chi_{\gamma \mu} (t_\lam) = \gamma^{\lam} \chi_{\mu}(t_\lam) \ \ \
\mu,\lam \in \Lambda^*_k,\ \gamma\in\Gamma
\end{equation}
where $\gamma^{\lam}$ is defined via the embedding $\Gamma \to T$.
We find,
$$ 
\sum_{\gamma \in \Gamma/\Gamma_\mu} 
\chi_{\gamma  \mu}(t_\lambda)^* =
\f{1}{\# \Gamma_\mu}
\chi_{\mu}(t_\lam)^*\ \sum_{\gamma \in \Gamma}\gamma^\lam . 
$$
The sum $\sum_{\gamma \in \Gamma}\gamma^\lambda$ is equal to 
$\#\Gamma$ if $\lambda$ is a weight of $T'=T/\Gamma$, and $0$ otherwise.  
Let $\Lambda_k^{'*} \subset \Lambda_k^*$ denote 
those level $k$ weights which are
weights for $T'$.  For any $\lambda \in \Lambda^{'*}_k$, the character
$\chi_\lambda$ descends to $T'$.  Therefore
$$ \chi(\M'(\Sig,\Co'_\mu)) = 
\f{\# T_{k+c}^{h-1}}{ \#
\Gamma_\mu \# \Gamma^{2h-1}} \sum_{ \lam \in {\Lambda}^{'*}_k} 
\eps(\phi,\lambda)\frac{ \# \Gamma_{\lambda}^{2h}}{
|J(t_\lambda)|^{2h-2}}\ \chi_{\mu}(t_\lambda)^* . $$
Finally, since 
$$ \chi_{\mu}(t_\lambda) = 
\chi_\mu(w_\gamma(t_\lam))=
\chi_{\mu}(\gamma t_\lambda),$$
and using that the factor $\eps(\phi,\lambda)\in\{0,1\}$ is invariant
under the action of $\Gamma$, we can re-write the result as a sum over
$\Lambda^{'*}_k/\Gamma$:
\begin{theorem}  Let $\mu \in \Lambda_k^*$ and  
$\M'(\Sig,\Co'_\mu)$ the moduli space of flat $G'$-connections on
$\Sig^h_1$, with holonomy around the boundary in the conjugacy class
$\Co'_\mu \subset G'$. The
 $\on{Spin_c}$-index of $\M'(\Sig,\Co'_\mu)$ is given by
$$
\chi(\M'(\Sig,\Co'_\mu)) 
= \f{ \# T_{k+c}^{h-1}}{\#
\Gamma_\mu \# \Gamma^{2h-2}} 
\sum_{ \lam \in {\Lambda}^{'*}_k/\Gamma}
\eps(\phi,\lambda)
\frac{ \# \Gamma_{\lambda}^{2h-1} }{
|J(t_\lambda)|^{2h-2}}\  \chi_{\mu}(t_\lambda)^*. $$
\end{theorem}

\begin{appendix}

\section{The cardinality of $T_l$}
\label{app:groups}
Suppose $G$ is simple and simply connected. For a proof 
of the following fact, see e.g. 
Beauville \cite[Remark 9.9]{be:co}, Bismut-Labourie 
\cite[Prop. 1.2, 1.3]{bi:sy}. 
\begin{lemma}\label{lem:bila} 
For any $l>0$ the number of elements in $T_l$ is given by 
$$ 
{\# T_l}=
l^{rank(G)} \Vol(T)^2=
l^{rank(G)}\,{\# Z(G)}\,\,\# (\ol{\mf{R}}/\ol{\mf{R}}_{long}).
$$
Here $\Vol(T)$ is the Riemannian volume of $T$ computed with respect
to the basic inner product, $\ol{\mf{R}}$ is the root lattice
and $\ol{\mf{R}}_{long}$ the sub-lattice generated by the long roots.
\end{lemma}
The table below gives the centers $Z(G)$, dual Coxeter numbers $c$, 
and the index of the long-root lattice $\ol{\mf{R}}_{long}$ inside 
the root lattice $\ol{\mf{R}}$ for all simple, simply connected compact 
Lie groups $G$. All of this information can be read off from 
tables in Bourbaki \cite{bo:li}. 
%
%
$$ \begin{array}{c||c|c|c|c|c|c|c|c|c} 
G& A_N & B_N & C_N & D_N & E_6 & E_7 & E_8 & F_4 & G_2 \\ 
\hline \hline\rule{0mm}{5mm}
c& N+1 & 2N-1 & N+1 & 2N-2 & 12 & 18 & 30 & 9 & 4 \\ \hline
\rule{0mm}{5mm}
& & & & \Z_4 \mbox{ ($N$ odd)}
& & & & & 
\\\rb{$Z(G)$}&
\rb{$\Z_{N+1}$}&
\rb{$\Z_2$}&
\rb{$\Z_2$}& \Z_2\times \Z_2 \mbox{ ($N$ even)}     
&\rb{$\Z_3$}&\rb{$\Z_2$}&\rb{0}&\rb{0}&\rb{0} 
\\
\hline
\rule{0mm}{6mm}
\#(\ol{\mf{R}}/\ol{\mf{R}}_{long})  
& 1 & 2 & 2^{N-1} & 1 & 1 & 1  & 1 & 4 & 3 
\\
\end{array}
$$

\section{Fusion of group valued Hamiltonian $G$-spaces}\label{sec:fusion}

In this appendix we collect some facts about fusion of Hamiltonian 
$G$-spaces with group-valued moment maps.
\begin{theorem}\cite[Theorem 6.1]{al:mom}
Let $G,H$ be compact Lie groups, and 
$(M,\om,(\Phi_1,\Phi_2,\Psi))$ a group valued Hamiltonian 
$G\times G\times H$-manifold. Let $M_{\on{fus}}$ be the same manifold 
with diagonal $G\times H$-action, $\Phi_{\on{fus}}=\Phi_1\Phi_2$, and 
$\om_{\on{fus}}=\om-\hh B(\Phi_1^*\theta,\Phi_2^*\olt)$. Then 
$(M_{\on{fus}},\om_{\on{fus}},(\Phi_{\on{fus}},\Psi))$ is a group valued 
Hamiltonian $G\times H$-manifold. 
\end{theorem}
The correction term $\hh B(\Phi_1^*\theta,\Phi_2^*\olt)$ will be 
loosely referred to as the ``fusion term''.   
If $M=M_1\times M_2$ is a product of
two $G\times H_i$-valued Hamiltonian spaces, we also write
$ M_1\fus M_2:=(M_1\times M_2)_{\on{fus}}$.  

Recall that if $G$ is a torus, a space with $G$-valued moment map
is just a symplectic manifold with a multi-valued moment map in the 
usual sense. Fusion of such spaces changes the symplectic form, 
but not the volume:
\begin{lemma}\label{lem:fusabelian}
Suppose $T$ is a torus, that $(M,\om,(\Phi_1,\Phi_2))$ a compact group
valued Hamiltonian $T\times T$-space, and
$(M_{\on{fus}},\om_{\on{fus}},(\Phi_{\on{fus}}))$ is the group valued
Hamiltonian $T$-space obtained by fusion. Then the symplectic volumes
of $M$ and $M_{\on{fus}}$ are the same.
\end{lemma}
\begin{proof}
This is a special case of a result for non-abelian groups proved in
\cite{al:du}. In the abelian case, the following much simpler argument
is available.  Notice that $M$ with diagonal $T$-action has moment map
$\Phi_{\on{fus}}=\Phi_1\Phi_2$ not only for the fusion form
$\om_{\on{fus}}$ but also for the original symplectic 2-form $\om$.
Suppose $t\in T$ is a (weakly) regular value of $\Phi_{\on{fus}}$, so that
$(\Phi_{\on{fus}})^{-1}(t)$
is a smooth sub-manifold and
$M_t=(\Phi_{\on{fus}})^{-1}(t)/T$ is an orbifold.  Since the pull-back
of the 2-form $\hh B(\Phi_1^*\theta,\Phi_2^*\theta)$ to 
$(\Phi_{\on{fus}})^{-1}(t)$
vanishes, the reduced symplectic forms are the same:
$\om_t=(\om_{\on{fus}})_t$. It follows that the two
Duistermaat-Heckman measures
$\m=\f{1}{n!}(\Phi_{\on{fus}})_*(|\om^n|)$ and
$\m_{\on{fus}}=\f{1}{n!}(\Phi_{\on{fus}})_*(|\om^n_{\on{fus}}|)$
coincide. Since the symplectic volume is the integral of the
Duistermaat-Heckman measure, the proof is complete.
\end{proof}
For any group-valued Hamiltonian $G$-space, the 2-form $\om$ is 
non-degenerate on the tangent space at any point in the identity 
level set. The following Lemma shows that fusion does not change 
the isotopy class of this symplectic structure. Its proof relies 
on the notion of {\em exponential} of a Hamiltonian space 
\cite{al:mom}: Let $\varpi\in \Om^2(\g)$ be the image of 
$\exp^*\eta$ under the homotopy operator $\Om^\star(\g)\to \
\Om^{\star-1}(\g)$. Then if $(M,\om_0,\Phi_0)$ is a Hamiltonian 
$G$-space in the usual sense, with $\Phi_0(M)$ contained in 
a sufficiently small neighborhood of $0$, then $(M,\om,\Phi)$ 
with $\om=\om_0+\Phi_0^*\varpi$ and $\Phi=\exp(\Phi_0)$ is 
a group valued Hamiltonian $G$-space. Conversely, if 
$(M,\om,\Phi)$ is a group-valued Hamiltonian $G$-space, 
any small neighborhood of $\Phi^{-1}(e)$ is obtained in this 
way. The 2-form $\varpi$ vanishes at $0$, so that 
$\om_m=(\om_0)_m$ for points in the zero level set. 
\begin{lemma}\label{lem:homotop}
Let $(M,\om,(\Phi_1,\Phi_2,\Psi))$ be a group valued Hamiltonian 
$G\times G\times H$-space, and $(M_{\on{fus}},\om_{\on{fus}},(\Phi_{\on{fus}},\Psi))$
its fusion. Let $m\in M$ be a point in the identity level set 
of $(\Phi_1,\Phi_2,\Psi)$. The symplectic 2-forms $\om|_m$ and 
$\om_{\on{fus}}|_m$ on $T_mM$ are isotopic through a path of symplectic
forms, invariant under the stabilizer group $(G\times H)_m$.
\end{lemma}

\begin{proof}
We may assume that $M$ is the exponential of a Hamiltonian 
$G\times G\times H$-space $(M,\om_0,(\Phi_{0,1},\Phi_{0,2},\Psi_0))$. 
Rescaling by $s>0$, we obtain a family of Hamiltonian spaces
$(M,s\om_0,(s\Phi_{0,1},s\Phi_{0,2},s\Psi_0))$, together with their 
exponentials. Let $\om_{\on{fus}}^s$ be the corresponding fusion forms. 
We claim that $s^{-1}\om_{\on{fus}}^s|_m$ give the required 
isotopy of symplectic forms. Indeed, each $\om_{\on{fus}}^s|_m$ is 
symplectic, and for $s\to 0$,  
$$
\om_{\on{fus}}^s|_m=s\om_0|_m-\hh B(\exp(s\Phi_{0,1})^*\theta,
\exp(s\Phi_{0,2})^*\olt)|_m=s\om_0|_m+O(s^2)
$$ 
showing that $\lim_{s\to 0} s^{-1}\om_{\on{fus}}^s|_m=\om_0|_m$.
\end{proof}

\section{Lifts of Weyl group elements} \label{app:lifts}
In this section we prove the following result about the 
embedding $\varphi:\ Z(G)\to W$ introduced earlier: 
\begin{proposition}\label{prop:lifts}
For any simply connected Lie group $G$ and any two elements
$\gamma_1,\gamma_2\in Z(G)$ of the center, the Weyl group
elements $w_j=\varphi(\gamma_j)\in W=N_G(T)/T$ lift to 
commuting elements $g_1,g_2\in N_G(T)$. 
\end{proposition}

\begin{proof}
Decomposing into irreducible factors we may assume $G$ is simple.  For
any simple group other than $G=D_N$ with $N$ even, the center $Z(G)$
is a cyclic group and the claim follows by choosing a lift of the
generator.  It remains to consider the case $G=D_N=\Spin(2N),\ N$ even
which has center $\Z_2\times\Z_2$.  We use the usual presentation
\cite{bo:li} of the root system of $D_N$ as the set of vectors $\pm
\eps_i \pm \eps_j,\ i\not=j$ in $\t^*\cong \R^N$, where $\eps_i$ are
the standard basis vectors of $\R^N$. The basic inner product on $D_N$
is the standard inner product on $\R^N$, and will be used to identify
$\t^*\cong\t$.  We choose simple roots $\alpha_j = \eps_j -
\eps_{j+1}, \ j = 1,\ldots,N-1$, and
$\alpha_N=\eps_{N-1}+\eps_N$. Then $\alpha_{\on{max}}=\eps_1+\eps_2$
is the highest root, and the fundamental alcove $\Alc$ is the polytope
defined by $\xi_1\ge \xi_2\ge\ldots\ge \xi_{N-1}\ge |\xi_N|$ and
$\xi_1+\xi_2\le 1$. The following four vertices of
$\Alc$ exponentiate to the central elements of $G$:
$$ 
\zeta_0=(0,\ldots,0),\ \zeta_1=(1,0,\ldots,0),\ \zeta_2=\hh(1,\ldots,1),\ 
\zeta_3=\hh(1,\ldots,1,-1).
$$
To describe the homomorphism \eqref{eq:imbed}, consider the exceptional
element 
$$\rho/c=\f{1}{2N-2}(N-1,N-2,N-3,\ldots,1,0) .$$
We have
$$ \rho/c+\zeta_1=\f{1}{2N-2}(3N-3,N-2,N-3,\ldots,1,0)=w_1\rho/c+\lambda_1,$$
where $\lambda_1=2\eps_1\in \ol{\mf{R}}=\Lambda$, and 
$w_1=\varphi(\exp(\zeta_1))\in W$  is the 
Weyl group element, 
$$ w_1(\xi_1,\ldots,\xi_N)=(-\xi_1,\xi_2,\ldots,\xi_{N-1},-\xi_N).$$ 
Similarly 
$$ \rho/c+\zeta_2=\f{1}{2N-2}(2N-2,\,2N-3,\ldots,\,N-2,\,N-1)
=w_2\rho/c+\lambda_2,$$
with $\lambda_2=\sum_{j=1}^N \eps_j\in\Lambda$ (since $N$ is even)
and $w_2=\varphi(\exp\zeta_2)$ given as 
$$ w_2(\xi_1,\ldots,\xi_N)=(-\xi_N,-\xi_{N-1},\ldots,-\xi_2,-\xi_1).$$
We construct lifts $g_1,g_2$ of $w_1,w_2$ in two stages. First 
we lift to commuting elements $g'_1,g'_2 \in \SO(2N)$ given in terms of their 
action on $\R^{2N}$ as follows, 
\beq
g'_1(x_1,\ldots,x_{2N}) &=& (x_1,-x_2,x_3,x_4,
\ldots,x_{2N-3},x_{2N-2},x_{2N-1},-x_{2N}),\\
g'_2(x_1,\ldots,x_{2N}) &=& (x_{2N-1},-x_{2N},
x_{2N-3},-x_{2N-2},\ldots,x_3,-x_4,-x_1,x_2) .
\eeq
We claim that any two lifts $g_1,g_2$ to $\Spin(2N)$ still commute.
The transformations $g_1',g_2'$ restrict to rotations on the 
$2N-2$-dimensional subspace 
$E=\{x|\ x_2=x_{2N}=0\}$ and on its orthogonal complement. The restrictions
of $g_j'$ to the subspace $E$ lift to commuting transformations 
of $\Spin(E)$ since $g_1'|E$ is the identity. 
The restrictions to $E^\perp$ lift to commuting transformations 
of $\Spin(E^\perp)$ since $E^\perp$ is two-dimensional and 
$\Spin(2)=\U(1)$ is abelian. 
\end{proof}

\end{appendix}

\end{document}